\documentclass[10pt]{article}
\usepackage{latexsym}
\usepackage{epsfig}
\input epsf
 
\begin{document}
\title{Anisotropic umbilic points and Hopf's Theorem for surfaces with constant anisotropic mean curvature}
\author{By M{\footnotesize IYUKI} K{\footnotesize OISO}
and B{\footnotesize ENNETT} P{\footnotesize ALMER}}
\date{}

\maketitle

\newtheorem{theorem}{Theorem}[section]
\newtheorem{cor}{Corollary}[section]
\newtheorem{prop}{Proposition}[section]
\newtheorem{lemma}{Lemma}[section]
\newtheorem{condition}{Condition}[section]
\newtheorem{example}{Example}[section]
\newtheorem{definition}{Definition}[section]
\newtheorem{remark}{Remark}[section]
\newtheorem{conjecture}{Conjecture}[section]
\newtheorem{claim}{Claim}[section]
\newcommand{\rf}[1]{\mbox{(\ref{#1})}}

\renewcommand{\thefootnote}{\fnsymbol{footnote}}


\begin{abstract}
We show that  for elliptic parametric functionals whose Wulff shape is smooth and has strictly positive curvature, any surface with constant anisotropic mean curvature which is a topological sphere is a rescaling of the Wulff shape.
\end{abstract}
\section{Introduction} Let $\gamma:S^2\rightarrow {\bf R}_+$ be a ``reasonable" positive function on the
two-dimensional unit sphere $S^2$.
For a smooth, oriented immersed surface $X:\Sigma \rightarrow {\bf R}^3$ with unit normal $\nu$, we define a functional by
\begin{equation}
\label{F}
{\cal F}[X]=\int_\Sigma \gamma(\nu)\:d\Sigma\:,\end{equation}
where $d\Sigma$ is the area element of $X$.
We will impose a {\it convexity condition} on the functional by requiring that the map
\begin{equation}
\label{chi}
{\tilde \chi}:S^2\rightarrow {\bf R}^3\:,\qquad \nu\mapsto D\gamma +\gamma \nu\:,\end{equation}
defines a smooth, convex surface $W:={\tilde \chi}(S^2)$. This surface is called the {\it Wulff shape}. Wulff's Theorem states that for all closed surfaces $S$ enclosing the same volume as $W$, ${\cal F}[W]\le {\cal F}[S]$ holds, so that $W$ solves the isoperimetric problem for this functional.
For example, if $|\cdot |$ is a smooth norm on ${\bf R}^3$ with dual
norm $|\cdot |_*$, then the functional obtained from using the density $\gamma(\nu):=
|\nu|$, satisfies the convexity condition and has the Wulff shape $W:=\{x\:\bigr| \:|x|_*=1\}$.

Now let $X_t=X+t \delta X+{\mathcal O}(t^2)$ be a smooth, compactly supported variation of $X$. The anisotropic mean curvature $\Lambda$ is defined by the first variation formula
\begin{equation}
\label{first}
\delta {\mathcal F}[X] :=\partial_t {\mathcal F}[X_t]_{t=0}=-\int_\Sigma \Lambda \delta X\cdot d\Sigma\:.\end{equation}
Since
$$\delta {\rm vol}[X]=\int_\Sigma \delta X\cdot d\Sigma\:,$$
the equation $\Lambda\equiv $ constant characterizes critical points of ${\cal F}$ with the enclosed volume constrained to be a constant.

A consequence of the convexity condition is that the equation for constant anisotropic mean curvature (CAMC) surfaces is absolutely elliptic in the sense of Hopf \cite{H}. In particular, the equation for prescribed anisotropic mean curvature possesses a Maximum Principle analogous to the well known one for CMC
(constnat mena curvature) surfaces. Since the Maximum Principle is one of the most important analytic tools for dealing with CMC surfaces, it is not surprising to see that many results for CMC surfaces have natural extensions to CAMC surfaces. The isoperimetric property of the Wulff shape is one such example.
%
Generalizing the Barbosa-do Carmo theorem, it was shown in \cite{P1997} that the only closed, stable CAMC surface is the, up to homothety, the Wulff shape. Also, generalizing the Alexandrov Theorem, it was recently shown in \cite{HL2007a} that the only closed, embedded CAMC surfaces are rescalings of the Wulff shape.

In this paper, we will show the following:
\begin{theorem} \label{main}
Assume the convexity condition holds for the functional ${\cal F}$. Then the only closed, genus zero surfaces with constant anisotropic mean curvature are rescalings of the Wulff shape $W$.
\end{theorem}

Of course, when $\gamma \equiv 1$, this gives Hopf's famous result that the only CMC topological spheres are round.  Recently, two interesting partial results for the anisotropic case have appeared. One, due to Giga and Zhai, \cite{GZ}, roughly states that the result holds for functionals which are sufficiently close in the $C^2$ topology, to the area functional. The other, due to He and Li, proves the result under the assumption that a second invariant besides the anisotropic mean curvature is also constant. This second invariant is Tracee$_\Sigma (d{\tilde \chi} \circ d\nu\circ J)$, where $J$ is the almost complex structure of the surface. The constancy of this invariant together with the constancy of $\Lambda$ is equivalent to the holomorphicity of a type of Hopf differential.
\section{Anisotropic Umbilic Points} Let $X:\Sigma \rightarrow {\bf R}^3$ be an oriented surface. At a point $p\in \Sigma$, we can consider the sphere $S_p$ which is in oriented contact with the surface at $p$ and which has the same mean curvature as the surface has at $p$. The sphere $S_p$ is called the {\it central sphere} of the immersion at $p$. If $S_p$ and the surface have contact of order at least two at $p$, then $p$ is an umbilic point of $X$.

We now consider a fixed Wulff shape $W$. Recall that since $W$ is convex, the unit normals to $W$ are in one to one correspondence with the points in $S^2$. At each point in $p \in \Sigma$ we can consider the surface $\omega_p$ which is the unique rescaling of $W$ that is in oriented contact with the surface at $p$ and has the same anisotropic mean curvature as the surface has at $p$.  We will call $p$ an {\it anisotropic umbilic}, (A-umbilic),  if $\omega_p$ and $X$ have at least second order contact at $p$.  We will show that either the surface is made up entirely of A-umbilics or these points are isolated and we can associate an integer to each of them.

Set $\chi={\tilde \chi}\circ \nu$.
$\chi: \Sigma \to W$ is called the anisotropic Gauss map of $X$.
A local expression for the anisotropic mean curvature $\Lambda$ is,
\begin{equation}
\label{AMC}
\Lambda :=-{\rm Trace}_{\Sigma}\:d\chi.
\end{equation}
 The condition that a point $p$ is an A-umbilic is that
\begin{equation}
\label{dif}\bigl(d\chi+(\Lambda/2)dX\bigr)_p=0\:.\end{equation}
We can assume, by making a translation in ${\bf R}^3$ if necessary, that
\begin{equation}
\label{norm}
\bigl(\chi+(\Lambda/2)X\bigr)_p=0.
\end{equation}
By \rf{dif}, we have that at any A-umbilic point, the Gaussian curvatures of $\Sigma$ and $W$ satisfy,  $K_\Sigma =(\Lambda^2/4)K_W>0$, so it follows that near $p$, the anisotropic Gauss map $\chi$ is a local diffeomorphism.
On the other hand, the Gauss map of $W$ is a global diffeomorphism and so near
$p\in \Sigma$ and near $\chi(p)\in W$, both surfaces can be parameterized over $S^2$ by the inverses of their Gauss maps. For $W$, the map ${\tilde \chi}$ given in \rf{chi} is exactly
this parameterization. If $q$ denotes the support function of $X$, then
$X=Dq+q\nu$ locally parameterizes the surface.
 From \rf{AMC}, the equation that the anisotropic mean curvature is constant is
expressed on $S^2$ as
\begin{equation}
\label{EL}
{\rm Trace}_{S^2} (D^2\gamma+\gamma I)(D^2q+qI)^{-1}=-\Lambda\equiv {\rm constant}\:.
\end{equation}
Since $K_\Sigma>0$ holds near $p$, the matrix $(D^2q+qI)^{-1}$ is positive definite near $p$ and \rf{EL} can be considered as a linear elliptic equation $E[\gamma]=-\Lambda.$
Clearly $E[(\Lambda/2)q]=\Lambda$ and so $w:=\gamma +(\Lambda /2)q$ satisfies
$E[w]=0$.

Note also that, from \rf{norm},  we have
$$(Dw+w\nu)_{\nu(p)}=0\:.$$
In particular, since $Dw$ and $\nu$ are perpendicular,
\begin{equation}
\label{norm2}
w(\nu(p))=0\:,\qquad Dw_{\nu(p)}=0\:.
\end{equation}

As in \cite{O}, we next introduce local coordinates near $\nu(p)$ in $S^2$ using central projection.
For $\nu$ near $\nu(p)$, let $y=\pi(\nu)$ be the intersection of the line through the origin
of ${\bf R}^3$ and $\nu$ with $T_{\nu(p)}S^2$.
For an orthogonal coordinate $(y_1, y_2)$ in $T_{\nu(p)}S^2$, let $\rho= \sqrt{y_1^2+y_2^2}$.
For a function $f$ on $S^2$, define ${\underline f}:=(1+\rho^2)^{1/2} f\circ \pi^{-1}$. Then,
there holds, (equation (4.1) of \cite{O}),
\begin{equation}
\label{f}
(1+\rho^2)^{1/2}({\underline f}_{y_i y_j})=\bigl(D^2f+fI \bigr)\:.\end{equation}
In these coordinates, the equation $E[w]=0$ has an expression
\begin{equation}
\label{eq1}a{\underline w}_{y_1 y_1} -2b{\underline w}_{y_1 y_2}+c{\underline w}_{y_2 y_2}=0\:,\end{equation}
for suitable functions $a,b$ and $c$.

As in \cite{H}, there is a linear change of coordinates $\xi_i =c_{i1}x_1+c_{i2}x_2$, with $(c_{ij})$ a constant matrix,
such that the previous partial differential equation takes the form
\begin{equation}
\label{eq2}a_1{\underline w}_{\xi_1\xi_1}-2b_1{\underline w}_{\xi_1\xi_2}+c_1{\underline w}_{\xi_2\xi_2}=0\:\end{equation}
with
\begin{equation}
\label{e3}a_1(0)=1=c_1(0)\:,\qquad b_1(0)=0\:. \end{equation}
By a theorem of Bers \cite{B},  there exists a homogeneous polynomial $P$
of degree $N$, $P$ not identically zero,  such that for all $\epsilon \in (0,1)$
\begin{equation}
\label{bigoi}
{\underline w}(\xi)=P(\xi)+{\cal O}(|\xi|^{N+\epsilon})\:,
\end{equation}
\begin{equation}
\label{bigoii}
{\underline w}_{\xi_i}=P_{\xi_i}+{\cal O}(|\xi|^{N-1+\epsilon})\:,\quad i=1, 2,
\end{equation}
\begin{equation}
\label{bigoiii}
{\underline w}_{\xi_i \xi_j}=P_{\xi_i \xi_j}+{\cal O}(|\xi|^{N-2+\epsilon})\:,\quad 1\le i,j \le 2,
\end{equation}
and
\begin{equation}
\label{Lap}
P_{\xi_1\xi_1}+P_{\xi_2\xi_2}=0\:
\end{equation}
holds on a neighborhood of $0$. Note that by \rf{norm2}, we have
$$0={\underline w}(0)={\underline w}_{\xi_i}(0),\:i=1,2\:.$$
It follows that $N\ge 2$ holds since $N\le 1$ together with \rf{bigoi} and \rf{bigoii} implies that $P\equiv 0$
by letting $\xi\rightarrow 0$.

We can also note that since $\xi=0$ corresponds to an A-umbilic, $N\ge 3$ holds since the left hand side of
\rf{bigoiii} vanishes when $\xi=0$. (If $N=2$ then for some $(i,j)$, $P_{\xi_i \xi_j}(0)\ne 0$ and letting
$\xi\rightarrow 0$ in \rf{bigoiii} gives a contradiction.)

Let $\zeta =\xi_1+\sqrt{-1} \xi_2$. By \rf{Lap}, $P_\zeta $ is a holomorphic function and so we can write
$P_{\zeta \zeta}=:\zeta ^{(N-2)}G(\zeta)$ where $G$ is a holomorphic function of $\zeta$ which is non vanishing in a
neighborhood of $\zeta=0$.
We obtain from \rf{bigoiii},
\begin{equation}
\label{bo}
{\underline w}_{\zeta \zeta}=\zeta^{N-2}\biggl[G(\zeta)+{\cal O}(|\zeta|^{\epsilon})\biggr]\:,\:\forall\epsilon \in (0,1)\:.
\end{equation}

Suppose there exists a sequence $\zeta_\mu \rightarrow 0$ with ${\underline w}_{\zeta \zeta}(\zeta_\mu)=0$, $\mu=1,2,3, \cdots$.
Then we obtain from \rf{bo},
$$G(\zeta_\mu)={\cal O}(|\zeta_\mu|^{\epsilon})\:,\epsilon \in (0,1),$$
which is a contradiction since $G(0)\ne0$. This shows that $\zeta=0$ is an isolated zero of the matrix $({\underline w}_{\xi_i \xi_j})$.
Using \rf{f} with $f=w$,  we see that  the A-umbilic at  $p$ is isolated.

\section{Indices}
In this section we show the following.
\begin{prop} Let $X:\Sigma \rightarrow {\bf R}^3$ be a CAMC surface which is not a rescaling of the Wulff shape. Let $p\in \Sigma$ be an A-umbilic and let $F$ be an eigendirection field for $D^2[\gamma +(\Lambda /2)q]+[\gamma +(\Lambda /2)q]I$ defined near $p$. Then the rotation index of $F$ around $p$ is negative.
\end{prop}
{\it Proof.$\:\:$}We will precisely describe the coordinate change in going from \rf{eq1} to \rf{eq2}.

Let
$${\cal L}_y= \left( \begin{array}{cc}
                       a(y) &  -b(y) \\
                     -b(y) & c(y)
                         \end{array} \right)\:.
$$
Let $\Lambda_i^2$, ($\Lambda_i>0$), $i=1,2$ be the eigenvalues of the  symmetric positive definite matrix
${\cal L}_0$. Then,  for some rotation matrix
$${\cal R}=\left( \begin{array}{cc}
                       \cos \vartheta  &- \sin \vartheta \\
                      \sin \vartheta &   \cos \vartheta
                         \end{array} \right)\:,$$
there holds ${\cal R}{\cal L}_0{\cal R}^{-1} =$ diagonal$(\Lambda_1^2, \Lambda_2^2)$.
It follows that if we make the coordinate transformation $t_1=(\cos \vartheta) y_1-(\sin \vartheta) y_2$, $t_2=(\sin \vartheta )y_1+(\cos \vartheta) y_2$, then the equation
$a(0)f_{y_1y_1}-2b(0)f_{y_1y_2}+c(0)f_{y_2y_2}=0$ is transformed into
$\Lambda_1^2 f_{t_1t_1}+\Lambda_2^2 f_{t_2t_2}=0$. Finally, the transformation
$\xi_i:=t_i/\Lambda_i$ changes this last equation into the Laplace equation
$f_{\xi_1\xi_1}+f_{\xi_2\xi_2}=0$.

The polynomial $P$ found above, therefore satisfies $\Lambda_1^2P_{t_1t_1}+
\Lambda_2^2P_{t_2t_2}=0$ and the Hessians $(P_{t_\alpha t_\beta})$ and
$(P_{y_i y_j})$ are related by
\begin{equation}
\label{hess}{\cal R}^{-1} (P_{t_\alpha t_\beta}){\cal R}=(P_{y_i y_j})\:.\end{equation}


From \rf{bigoiii}, we get
\begin{equation}
\label{bigoiv}
{\underline w}_{y_i y_j}=P_{y_i y_j}+{\cal O}(|y|^{N-2+\epsilon})\:,\quad 1\le i,j \le 2\:.
\end{equation}
By \rf{hess} the rotation index of the eigendirections of $(P_{t_\alpha t_\beta})$
and $(P_{y_i y_j})$ are the same since the eigendirections for one of the matrices differs from the eigendirection for the other by a fixed rotation. Therefore it is enough to show
that the rotation index of the eigendirections of $(P_{t_\alpha t_\beta})$
are negative.

We write $\xi_1+i\xi_2=\rho e^{i\theta}$, then a calculation shows
$$(P_{t_\alpha t_\beta})=N(N-1)\rho^{N-2}\left( \begin{array}{cc}
                       \frac{ \cos (N-2)\theta}{\Lambda_1^2} & \frac{-\sin(N-2)\theta}{\Lambda_1\Lambda_2} \\
                     \frac{ -\sin (N-2)\theta}{\Lambda_1\Lambda_2} &\frac{ -\cos (N-2)\theta}{\Lambda_2^2}
                         \end{array} \right)\:.$$
   Let
\begin{eqnarray*}
\Delta:&=&\{(\cos^2(N-2)\theta)\bigl(\frac{1}{\Lambda_2^2}-\frac{1}{\Lambda_1^2}\bigr)^2
   +\frac{4}{\Lambda_1^2\Lambda_2^2}\}^{1/2}\\
   &=&\{(\cos^2(N-2)\theta)\bigl(\frac{1}{\Lambda_1^2}+\frac{1}{\Lambda_2^2}\bigr)^2+\frac{4\sin^2(N-2)\theta}{\Lambda_1^2\Lambda_2^2}\}^{1/2}\\
   &\ge& \frac{2}{|\Lambda_1\Lambda_2|}\qquad (*)\:.
   \end{eqnarray*}
Then the eigenvalues are given by
\begin{equation}
\label{lambda}
\lambda_{\pm}=\frac{N(N-1)}{2}\rho^{N-2}\biggl[\bigl(\cos(N-2)\theta)\bigl(\frac{1}{\Lambda_1^2}-\frac{1}{\Lambda_2^2}\bigr)\pm \Delta \biggr]  \:.\end{equation}
If $(x,y)$ is an eigenvector belonging to $\lambda_+$, we easily obtain
\begin{eqnarray*}
\tan \Phi:=\frac{y}{x}&=& \frac{\Lambda_1\Lambda_2}{2}\biggl(\cot((N-2)\theta)\bigl(\frac{1}{\Lambda_1^2}+\frac{1}{\Lambda_2^2}\bigr)  -\frac{\Delta}{\sin ((N-2)\theta)}\biggr)\\
&=&\frac{1}{2}\biggl(\cot((N-2)\theta) \bigl(\frac{\Lambda_2}{\Lambda_1}+\frac{\Lambda_1}{\Lambda_2}\bigr)
-\frac{\{(\cos^2(N-2)\theta)\bigl(\frac{\Lambda_2}{\Lambda_1}-\frac{\Lambda_1}{\Lambda_2}\bigr)^2
   +4\}^{1/2}}{\sin((N-2)\theta)}\biggr)\:.
   \end{eqnarray*}

It follows that the winding number of the eigendirection corresponding to $\lambda_+$, is
$$\frac{1}{2\pi}\int_0^{2\pi}
d\biggl(\arctan\biggl(\frac{1}{2}\biggl(\cot((N-2)\theta) \bigl(\frac{\Lambda_2}{\Lambda_1}+\frac{\Lambda_1}{\Lambda_2}\bigr)
-\frac{\{(\cos^2(N-2)\theta)\bigl(\frac{\Lambda_2}{\Lambda_1}-\frac{\Lambda_1}{\Lambda_2}\bigr)^2
   +4\}^{1/2}}{\sin((N-2)\theta}\biggr)\biggr)\biggr)=-\frac{N-2}{2}\:.$$
  (This is shown in the Appendix.)
Recall that we have shown above that $N\ge 3$ holds.
The right hand side of the above equality is ngative.

We next show that the rotation index of an eigendirection of $({\underline w}_{y_i y_j})$ at an A-umbilic is equal to the index of an eigendirection of $(P_{y_i y_j})$ at the same point.

 Write $\Omega=({\underline w}_{y_i y_j})$, ${\cal P}:=(P_{y_i y_j})$. The eigenvalues of ${\cal P}$ are given in \rf{lambda} and the corresponding eigenvectors are orthogonal since ${\cal P}$ is self-adjoint.  Set $\lambda_+=\lambda_1$, $\lambda_-=\lambda_2$. Suppose ${\cal P}E_1= \lambda_1E_1$,  ${\cal P}E_2=\lambda_2E_2$ with $E_i\cdot E_j=\delta_{ij}$.
Let $V$ be a unit eigenvector of $\Omega$ with eigenvalue $\lambda =\lambda(\rho,\theta)$.
Write $V=:(\cos \alpha)E_1+(\sin \alpha )E_2$.  Then
\begin{eqnarray*}\lambda \bigl((\cos \alpha)E_1+(\sin \alpha )E_2\bigr)
&=&\Omega \bigl((\cos \alpha)E_1+(\sin \alpha )E_2\bigr)\\
&=&\bigl({\cal P}+{\cal O}(\rho^{N-2+\epsilon})\bigr) \bigl((\cos \alpha)E_1+(\sin \alpha )E_2\bigr)\\
&=&\lambda_1(\cos \alpha)E_1+\lambda_2(\sin \alpha )E_2\\
&&+{\cal O}(\rho^{N-2+\epsilon})\bigl((\cos \alpha)E_1+(\sin \alpha )E_2\bigr)\:.
\end{eqnarray*}
This implies, using \rf{lambda}, that

\begin{equation}
\label{e1}
\cos \alpha \bigl(\frac{\lambda}{\rho^{N-2}}\bigr)
=\cos \alpha\biggl(\frac{N(N-1)}{2}[\cos((N-2)\theta)\bigl(\frac{1}{\Lambda_1^2}-\frac{1}{\Lambda_2^2}\bigr)+\Delta]+{\cal O}(\rho^\epsilon)\biggr)\:,\end{equation}

\begin{equation}
\label{e2}\sin \alpha \bigl(\frac{\lambda}{\rho^{N-2}}\bigr)
=\sin \alpha\biggl(\frac{N(N-1)}{2}[\cos((N-2)\theta)\bigl(\frac{1}{\Lambda_1^2}-\frac{1}{\Lambda_2^2}\bigr)-\Delta]+{\cal O}(\rho^\epsilon)\biggr) \:,\end{equation}
where $\Lambda_i$ are constants. It follows that either $\cos \alpha \rightarrow 0$ or $\sin \alpha \rightarrow 0$
as $\rho\rightarrow 0$. To see this, assume that neither of these limits hold. Since the circle is compact, we can find a sequence of points
$(\rho_i, \theta_i)$ with $\rho_i\rightarrow 0$ and $\cos(\alpha(\rho_i,\theta_i))\rightarrow A\ne 0$, $\sin(\alpha(\rho_i,\theta_i))\rightarrow B\ne 0$. We can use this to cancel the factors $\cos \alpha$ and $\sin \alpha$ in  \rf{e1} and \rf{e2} to obtain:
$$\frac{\lambda(\rho_i, \theta_i)}{\rho_i^{N-2}}- \frac{N(N-1)}{2}[\cos((N-2)\theta_i)\bigl(\frac{1}{\Lambda_1^2}-\frac{1}{\Lambda_2^2}\bigr)+\Delta (\theta_i)]={\cal O}(\rho_i^\epsilon)$$
and
$$\frac{\lambda(\rho_i, \theta_i)}{\rho_i^{N-2}}- \frac{N(N-1)}{2}[\cos((N-2)\theta_i)\bigl(\frac{1}{\Lambda_1^2}-\frac{1}{\Lambda_2^2}\bigr)-\Delta (\theta_i)]={\cal O}(\rho_i^\epsilon).$$
Subtracting, we obtain
$$\Delta (\theta_i)={\cal O}(\rho_i^\epsilon)\:,$$
which is impossible because of ($\ast$).

  We will assume $\sin \alpha\rightarrow 0$,  the other case is similar.

It follows that  $V$ and $E_1$ are asymptotically parallel as $\rho \rightarrow 0$ and hence $V$ and $E_1$ have the same winding numbers about $\rho=0$. To see this,
write $E_1=(\cos \mu(\rho, \theta), \sin  \mu(\rho, \theta))$, $E_2=(-\sin  \mu(\rho, \theta),
\cos  \mu(\rho, \theta))$, then $V=(\cos (\alpha+\mu), \sin(\alpha+\mu))$. The rotation index of $V$ is the integer $J$ given by
$$2\pi J=\lim_{\rho \rightarrow 0} \oint (d\alpha +d\mu)\:.$$
Choose $\rho_0 \approx 0$ such that  $0<\rho <\rho_0$ implies $|\sin \alpha|<1/2$.
For some $\rho$, $0<\rho <\rho_0$, if $\alpha(\rho,0)$ is continued along the circle of
radius $\rho$, then we arrive at $\alpha(\rho,2\pi)= \alpha(\rho,0)+\pi m$ for some integer $m$. However, $m=0$ must hold, otherwise there is a value between $\alpha(\rho,0)$
and $\alpha(\rho,0)+\pi m$ where $\sin \alpha=1$ holds, giving a contradiction.
Therefore,
$$\lim_{\rho\rightarrow 0}\oint d\alpha \rightarrow 0\:.$$
It follows that the index of an  eigendirection of $({\underline w}_{y_i y_j})$ and the index of an eigendirection of $(P_{y_i  y_j})$ are the same at $\xi=y=0$ and so both are negative. {\bf q.e.d.}\\[2mm]

\section{Proof of Main Result}

We will now assume that $\Sigma$ is a closed, genus zero surface with constant anisotropic mean curvature $\Lambda$. The reader can easily verify that the product of a positive definite  symmetric matrix, in this case  $D^2\gamma+\gamma I$, and another symmetric matrix, in this case $d\nu$, has real eigenvalues. By considering the discriminant of the characteristic polynomial of $(D^2\gamma+\gamma I)\cdot d\nu$, one can then see that  $(\Lambda^2/4)-K_\Sigma/K_W \ge 0$ holds. So $\Lambda\equiv 0$ would imply that $\Sigma$ has non positive curvature which is impossible for a closed surface in 3-space.

In order to prove the theorem, we will apply a well known result on the sum of the indices of a direction field, \cite{H},  to an eigendirection field of $D^2[\gamma +(\Lambda /2)q]+[\gamma +(\Lambda /2)q]I$.  Although $D^2 \gamma+\gamma I$ is globally well defined, the endomorphism field $D^2q+qI$ is
undefined at points where the curvature $K_\Sigma$ vanishes.
\begin{lemma}
Let $C$ be the closed subset of $\Sigma$ where $K_\Sigma=0$ holds. Let $v$ be an eigendirection field of $D^2[\gamma +(\Lambda /2)q]+[\gamma +(\Lambda /2)q]I$. Then $v$ can be continued continuously across $C$.
\end{lemma}
{\it Proof.$\:\:$} We first show that there is a neighborhood $C'$ of $C$ which does not contain any umbilic points or A-umbilic points. If $C$ contained an umbilic point, then the principal curvatures at that point would satisfy $k_1-k_2=0=k_1k_1$ so $d\nu$ would vanish. This would make
\rf{AMC} with $\Lambda \ne 0$ impossible. It follows that there is a neighborhood $C_1$ of $C$ where which is free of umbilics.  At any A-umbilic,
we have $(\Lambda^2/4)-K_\Sigma/K_W=0$. Therefore $C_2:=\{(\Lambda^2/4)-K_\Sigma/K_W > \Lambda^2/8\}$, gives a neighborhood of $C$ containing no A-umbilics. We let $C'=C_1\cap C_2$.

We will show that the endomorphism field $K(D^2[\gamma +(\Lambda /2)q]+[\gamma +(\Lambda /2)q]I)$ extends continuously to $C$ and that this endomorphism field has no singularities in $C$ (singularites here includes the possiblitity that the field vanishes). Since $D^2[\gamma +(\Lambda /2)q]+[\gamma +(\Lambda /2)q]I$ and $K(D^2[\gamma +(\Lambda /2)q]+[\gamma +(\Lambda /2)q]I)$ have the same eigengirection fields on $C'\setminus C$, we can extend the eigendirection fields
of $D^2[\gamma +(\Lambda /2)q]+[\gamma +(\Lambda /2)q]I$ to $C$ by using the eigendirection fields of $K(D^2[\gamma +(\Lambda /2)q]+[\gamma +(\Lambda /2)q]I)$.

We work at a point in $x\in C'\setminus C$. At $x$, choose an orthonormal frame consisting of principal directions. With respect to this frame we write
$$ D^2\gamma+\gamma I=\left( \begin{array}{cc}
                       a_{11}  &a_{12}\\
                      a_{12} &   a_{22}
                         \end{array} \right)\:,\:\: D^2q+q I=\left( \begin{array}{cc}
                       -1/k_1  &0\\
                      0 &   -1/k_2
                         \end{array} \right).$$
Straightforward calculations then gives
\begin{equation}
\label{adn}
(D^2\gamma+\gamma I)\cdot d\nu=\left( \begin{array}{cc}
                       -k_1a_{11}  &-k_2a_{12}\\
                      -k_1a_{12} &   -k_2a_{22}
                         \end{array}\right)
                         \:,
\end{equation}
$$\Lambda=k_1a_{11}+k_2a_{22}\:,$$
and
\begin{equation}
 K(D^2[\gamma +(\Lambda /2)q]+[\gamma +(\Lambda /2)q]I)
=\left( \begin{array}{cc}
                    k_2(k_1a_{11}-k_2a_{22})/2 &Ka_{12}\\
                    Ka_{12} &k_1(k_2a_{22}-k_1a_{11})/2
                         \end{array}\right)\:.
\end{equation}
From \rf{adn}, we see that an A-umbilic on a  surface with constant non zero $\Lambda$ corresponds to a point where $k_1a_{11}-k_2a_{22}=0$ and $a_{12}=0$.  Since there are no umbilics in $C'$, the frame of principal directions is well defined on $C'$. Since both principal curvatures cannot simultaneously vanish and since there are no A-umbilics in $C'$, the previous matrix cannot vanish anywhere on $C'$.  and the expression above for $K(D^2[\gamma +(\Lambda /2)q]+[\gamma +(\Lambda /2)q]I)$ extends over $C$ and has well defined eigen-direction fields there. {\bf q.e.d.}

\vskip0.5truecm

\noindent {\it Proof of Theorem \ref{main}.} Assume that the surface is a topological sphere with constant anisotropic mean curvature which is not a rescaling of the Wulff shape. Then the anisotropic umbilic points are isolated.  We consider a direction field $F$ on $\Sigma$ which is given on $\Sigma \setminus C'$ as an eigendirection field of $D^2[\gamma +(\Lambda /2)q]+[\gamma +(\Lambda /2)q]I$ and is given on $C'$ as an eigendirection field of $K(D^2[\gamma +(\Lambda /2)q]+[\gamma +(\Lambda /2)q]I)$ as defined above. Then the rotation indices of the singularities of $F$ are all negative. However, the sum of the indices of a line field on a topological sphere is positive, which gives a contradiction.
{\bf q.e.d.}

\section{Acknowledgements}

The authors wish  to thank  Professor  Naoya
Ando who supplied us with some of the ideas used in Section 3.

The first author is partially supported by
Grant-in-Aid for Scientific Research (C) No. 19540217
of the Japan Society for the Promotion of Science.
The   second author  was partially funded by Fellowship S-08154 from the Japan Society for the Promotion of Science.  This work was also partially supported by Fundaci\'{o}n S\'{e}neca
project 04540/GERM/06, Spain. This research is a result of the activity
developed
within the framework of the Programme in Support of Excellence Groups of the
Regi\'{o}n de Murcia, Spain, by Fundaci\'{o}n S\'{e}neca, Regional
Agency for
Science and Technology (Regional Plan for Science and Technology 2007-2010).
\section{Appendix}
We show
$$\frac{1}{2\pi}\int_0^{2\pi}
d\biggl(\arctan\biggl(\frac{1}{2}\biggl(\cot((N-2)\theta) \bigl(\frac{\Lambda_2}{\Lambda_1}+\frac{\Lambda_1}{\Lambda_2}\bigr)
-\frac{\{(\cos^2(N-2)\theta)\bigl(\frac{\Lambda_2}{\Lambda_1}-\frac{\Lambda_1}{\Lambda_2}\bigr)^2
   +4\}^{1/2}}{\sin((N-2)\theta}\biggr)\biggr)\biggr)=-\frac{N-2}{2}\:.$$

   Let $z:=\Lambda_2/\Lambda_1$, $m:=N-2$,  $\psi:=m\theta$, then this formula is equivalent to
\begin{equation}
\label{f1}
-1/2= \frac{1}{2\pi} \int_0^{2\pi} \partial_\psi  \arctan \biggl(\frac{1}{2}\bigl([z+1/z]\cot \psi -\frac{\{4+[z-1/z]^2\cos^2\psi\}^{1/2}}{\sin \psi}\bigr)\biggr)\:d\psi\:.
\end{equation}
Note that the integrand has singulaities at each half odd integer multiple of $\pi$. In any interval free of half odd integer multiples of $\pi$,
the integrand can be computed as
$$
\partial_\psi  \arctan \biggl(\frac{1}{2}\bigl([z+1/z]\cot \psi -\frac{\{4+[z-1/z]^2\cos^2\psi\}^{1/2}}{\sin \psi}\bigr)\biggr)
=\frac{-z(z^2+1)}{(z^2-1)^2\cos^2\psi+4z^2}\:.$$
A standard table of integrals gives
$$\int\frac{d\psi}{p^2+q^2\cos^2\psi}=\frac{1}{p\sqrt{p^2+q^2}} \arctan \bigl(\frac{p\tan \psi}{\sqrt{p^2+q^2}}\bigr)\:.$$
Using this with $p^2=4z^2$ and $q^2=(z^2-1)^2$, we obtain
$$ \int \partial_\psi  \arctan \biggl(\frac{1}{2}\bigl([z+1/z]\cot \psi -\frac{\{4+[z-1/z]^2\cos^2\psi\}^{1/2}}{\sin \psi}\bigr)\biggr)\:d\psi
=-\frac{1}{2} \arctan \frac{2z\tan\psi}{z^2+1}\:.$$
Evaluating the antiderivative over the endpoints of the successive intervals $(0,\pi/2)$, $(\pi/2, 3\pi/2)$ and $(3\pi/2, 2\pi)$, gives
$$ \frac{1}{2\pi}\int^{2\pi}_0 \partial_\psi  \arctan \biggl(\frac{1}{2}\bigl([z+1/z]\cot \psi -\frac{\{4+[z-1/z]^2\cos^2\psi\}^{1/2}}{\sin \psi}\bigr)\biggr)\:d\psi
 =\frac{1}{2\pi}(-\pi)=-1/2\:,$$
 which proves \rf{f1}.


\begin{flushleft}
Miyuki K{\footnotesize OISO} \\
Department of Mathematics\\
Nara Women's University \& PRESTO, JST\\
Kita-Uoya Nishimachi\\
Nara 630-8506\\
Japan\\
E-mail: koiso@cc.nara-wu.ac.jp
\end{flushleft}

\begin{flushleft}
Bennett P{\footnotesize ALMER}\\
Department of Mathematics\\
Idaho State University\\
Pocatello, ID 83209\\
U.S.A.\\
E-mail: palmbenn@isu.edu
\end{flushleft}
\end{document}